\documentclass{amsart}
\usepackage{amsmath}
\usepackage{amssymb}
\usepackage{amsfonts}
\usepackage[dvips]{graphicx}

\newtheorem{claim}{Claim}
\newtheorem{theorem}[claim]{Theorem}

\newtheorem{observation}[claim]{Observation}

\theoremstyle{definition}

\begin{document}

\title{Dyck paths with coloured ascents}
\author{Andrei Asinowski$^\dag$ \and Toufik Mansour$^\ddag$}
\maketitle
\begin{center}
{$^\dag$Caesarea Rothschild Institute, University of Haifa, Haifa 31905,
Israel, +972-4-8288343}\\[2pt]
{$^\ddag$Department of Mathematics, University of Haifa, Haifa 31905,
Israel, +972-4-8240705}\\[6pt]

\verb1 andrei@cri.haifa.ac.il, toufik@math.haifa.ac.il1
\end{center}

\section*{Abstract} We introduce a notion of Dyck paths with coloured
ascents. For several ways of colouring, we establish bijections
between sets of such paths and other combinatorial structures, such
as non-crossing trees, dissections of a convex polygon, etc. In some
cases enumeration gives new expression for sequences enumerating
these structures.

\vspace{3mm}

\noindent \textsc{Keywords.} Dyck paths, non-crossing graphs,
dissections of polygon by diagonals.

\noindent \textsc{2000 Mathematics Subject Classification}: 05A05,
05A15.

\section{Introduction}
\subsection{Coloured Dyck paths}
\label{sec:introdyck} A {\it Dyck path of length $2n$} is a sequence
$P$ of letters $U$ and $D$, such that $\#(U)=\#(D)=n$ (where $\#$
means ``number of'') in $P$, and $\#(U) \geq \#(D)$ in any initial
subsequence of $P$. A Dyck path of length $2n$ is usually
represented graphically as a lattice path from the point $(0,0)$ to
the point $(2n, 0)$ that does not pass below the $x$-axis, where $U$
is the ``up step'' $(1, 1)$ and $D$ is the ``down step'' $(1, -1)$.
The set of all Dyck paths of length $2n$ will be denoted by
$\mathcal{D}(n)$. We shall also denote $\mathcal{D} = \bigcup_{n\geq
0}\mathcal{D}_n$. It is well known that $\vert \mathcal{D}(n) \vert$
is equal to the $n$-th Catalan number $C_n =
\frac{1}{n+1}\binom{2n}{n}$ (see \cite[Page 222, Exercise
19(i)]{S}).

A maximal subsequence of $k$ consecutive $U$'s
(that is, not preceded or followed by another $U$) in a Dyck path
will be called a {\it $k$-ascent}
and denoted by $U^k$.
Similarly, a maximal subsequence of $k$ consecutive $D$'s
will be
denoted by $D^k$.
The Dyck path $U^kD^k$
will be called {\it pyramid of length $2k$} and denoted by $\Lambda^k$.

In this paper we present a new generalization of Dyck paths.
Let $\mathfrak{L} = \{\mathcal{L}_0, \mathcal{L}_1, \mathcal{L}_2, \dots\}$
be a sequence of sets,
and let $a_k = \vert \mathcal{L}_k \vert$.
We colour all ascents in a Dyck path, when the set of colours for each $k$-ascent
is $\mathcal{L}_k$.
In this way we obtain {\it Dyck paths with ascents coloured by $\mathfrak{L}$}
(shortly {\it Dyck paths coloured by $\mathfrak{L}$}, or
{\it coloured Dyck paths}).
Each Dyck path $P$ produces thus $\prod a_i$ coloured Dyck paths,
when the product is taken over the lengths of all ascents in $P$.
Coloured Dyck paths will be denoted by capital letters with ``hat'', e.\ g.\ $\hat{P}$.
The pyramid $\Lambda^k$ with $U^k$ coloured by a specified colour $C \in \mathcal{L}(k)$
will be denoted by $\Lambda^k\langle C\rangle$.

The set of all Dyck paths of length $2n$
coloured by members of $\mathfrak{L}$
will be denoted by $\mathcal{D}^{\mathfrak{L}}(n)$.
We shall also denote
$\mathcal{D}^{\mathfrak{L}} = \bigcup_{n\geq 0}\mathcal{D}^{\mathfrak{L}}(n)$.
In order to obtain a general expression enumerating
$\vert\mathcal{D}^{\mathfrak{L}}(n)\vert$,
we note that any Dyck path can be presented uniquely in the form
\[ U^k D P_k D P_{k-1} \dots D P_1, \]
where $P_k, P_{k-1}, \dots, P_1$ are (possibly empty) Dyck paths.
Therefore $M(x)$, the generating function for the sequence $\{\vert
\mathcal{D}^{\mathfrak{L}}(n)\vert\}_{n\geq0}$, satisfies
\[ M(x)=a_0 + a_1 x M(x) + a_2 x^2 M^2(x) + \dots, \]
or
\begin{equation}
\label{eq:gen2}
M= A(xM),
\end{equation}
where $A(x) = \sum_{i\geq0} a_ix^i$ is the generating function for
the sequence $\{\vert \mathcal{L}_i \vert\}_{i\geq 0}$.

However, we shall consider rather $\mathcal{L}_k$'s than merely
$a_k$'s, and we shall establish bijections between
$\mathcal{D}^{\mathfrak{L}}(n)$ for some specific $\mathfrak{L}$,
and other combinatorial structures. For example, in our main result
$\mathcal{L}_k=\mathcal{D}(k)$, and thus $a_k=C_k =
\frac{1}{k+1}\binom{2k}{k}$. We shall show the bijection between
Dyck paths coloured in this way and non-crossing trees to be defined
in Section~\ref{sec:trees}.

\subsection{Non-crossing trees}
\label{sec:trees}

A {\it non-crossing tree} (an ``NC-tree'') {\it on $[n]$} is a
labeled tree which can be represented by a drawing in which the
vertices are points on a circle, labeled by $\{1, 2, \dots, n\}$
clockwise, and the edges are non-crossing straight segments.
The vertex $1$ will be also called {\it the root},
and we shall depict it as a top vertex.
Non-crossing trees have been studied by Chen et al. \cite{chen},
Deutsch et al. \cite{deutschFN, noy02}, Flajolet et al.
\cite{flajolet}, Hough \cite{hough}, Noy et al. \cite{N}, and
Panholzer et al. \cite{prodinger}. Denote the set of all NC-trees
on $[n]$ by $\mathcal{NC}(n)$.
It is well known that
$\vert \mathcal{NC}(n+1) \vert = \frac{1}{2n+1}\binom{3n}{n}$.

We shall use the following notion.
If $\{a, b\}$ is an edge of
$T$, we shall rather denote it by $(a, b)$,
always assuming that $a<b$.
For a vertex $v$, we denote
$N^-(v)=\{a\in V(T): (a, v)\in E(T) \}$ (the {\it in-edges incident to $v$}),
$N^+(v)=\{b\in V(T): (v, b)\in E(T) \}$ (the {\it out-edges incident to $v$});
$d^-(v)=\vert N^-(v) \vert$ (the {\it in-degree of $v$}),
$d^+(v)=\vert N^+(v) \vert$ (the {\it out-degree of $v$}).

Consider the NC-trees on $[n]$ with the property: For each
vertex $v \not = 1$, we have $d^-(v) = 1$.
Such trees will be called {\it non-crossing out-trees} (``NCO-trees'');
the set of NCO-trees on $[n]$ will be denoted by $\mathcal{NCO}(n)$.

We also denote $\mathcal{NC} = \bigcup_{n\geq 0}\mathcal{NC}(n+1)$ and
$\mathcal{NCO} = \bigcup_{n\geq 0}\mathcal{NCO}(n+1)$.

\subsection{The results}
\label{seq:results}
We establish bijections between $\mathcal{D}^{\mathfrak{L}}(n)$
and other combinatorial structures, for a few specific
choices of $\mathfrak{L}$.
The main result is the following theorem:
\begin{theorem}
\label{the:dyckdyck}
There is a bijection between
the set of Dyck paths of length $2n$
with $k$-ascents coloured by
Dyck paths of length $2k$
and the set of
non-crossing trees on $[n+1]$.
\end{theorem}

Other results are special cases and variations of this theorem.
Substituting the generation function of $\mathfrak{L}$
in \eqref{eq:gen2}
enables us to
enumerate easily
the combinatorial structures
being in bijection with $\mathcal{D}^{\mathfrak{L}}(n)$.

\section{Dyck paths coloured by Dyck paths}
\label{sec:dyck}
Let $\mathfrak{D} = \{\mathcal{D}(0), \mathcal{D}(1), \mathcal{D}(2), \dots \}$.
In this Section we consider $\mathcal{D}^{\mathfrak{D}}(n)$ --
the set of Dyck paths
of length $2n$
with $k$-ascents coloured by Dyck paths of length $2k$,
i.\ e.\ we take $\mathcal{L}_k = \mathcal{D}(k)$.

\begin{figure}[ht]
$$\resizebox{60mm}{!}{\includegraphics{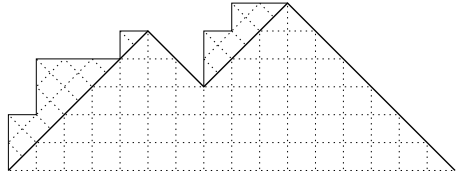}}$$
\caption{A Dyck path with $k$-ascents coloured by Dyck paths of length $2k$.}
\label{fig:dyckdyck}
\end{figure}

First we introduce a convenient way to depict thus coloured Dyck paths.
Given a $k$-ascent $U^k$ coloured by a Dyck path $C$ of length $2k$,
we draw a copy of $C$, rotated by $45^\circ$ and scaled by $1/\sqrt{2}$,
between the endpoints of $U^k$.
Figure~\ref{fig:dyckdyck} presents in this way
the Dyck path $U^5D^2U^3D^6$ with $U^5$ coloured by $UUDUUDDDUD$
and $U^3$ coloured by $UUDUDD$.

\subsection{Enumeration}
Let us enumerate $\mathcal{D}^{\mathfrak{D}}(n)$.
The generating function for $\{ \vert\mathcal{D}(n)\vert\}_{n \geq 0}$ is
\[ \frac{1-\sqrt{1-4x}}{2x} = 1 + x + 2x^2 + 5 x^3 + \dots . \]
Substituting this in \eqref{eq:gen2}, we obtain $$M =
\frac{1-\sqrt{1-4xM}}{2xM}.$$ After simplifications, we have $M-1 =
x M^3$. Denoting $L = M-1$ and applying Lagrange's inversion formula
(see~\cite[Section 5.4]{S} and \cite[Section 5.1]{wilf}) on $L = x
(L+1)^3$, we get that the coefficient of $x^n$ in $L$ is
\[ [x^n]L = \frac{1}{n}[L^{n-1}](L+1)^{3n}=
\frac{1}{n}\binom{3n}{n-1}=\frac{1}{2n+1}\binom{3n}{n}.
\] Thus we have $\vert \mathcal{D}^{\mathfrak{D}}(n) \vert = \vert
\mathcal{NC}(n+1) \vert$.

In this Section will shall construct,
for each $n\geq 0$, a bijective function
$\varphi_n: \mathcal{D}^{\mathfrak{D}}(n) \rightarrow \mathcal{NC}(n+1)$.
It will be presented as a restriction
of a bijective function $\varphi: \mathcal{D}^{\mathfrak{D}}\rightarrow \mathcal{NC}$.
The function $\varphi$ will be constructed by the following steps:
In Subsection~\ref{sec:decpath}
we describe a recursive procedure of decomposing a Dyck path into pyramids.
In Subsection~\ref{sec:dycknco} we construct
a bijection $\vartheta: \mathcal{D} \rightarrow \mathcal{NCO}$.
In Subsection~\ref{sec:bij} we first define $\varphi$
for coloured pyramids and then, using observations from Subsection~\ref{sec:decpath},
for all Dyck paths coloured by $\mathfrak{D}$.
All by all, this will give us the function $\varphi$,
and we shall also show that it is bijective.

\subsection{Decomposition of a Dyck path into pyramids}
\label{sec:decpath}
Let $P$ be a Dyck path.
Recall that it can be presented uniquely in the form
\begin{equation}
\label{eq:dec}
P = U^k D P_k D P_{k-1} \dots D P_1,
\end{equation}
where $P_k, P_{k-1}, \dots, P_1$ are (possibly empty) Dyck paths.
We say that $\Lambda^k$ is {\it the base pyramid of $P$} and that
$P_k, P_{k-1}, \dots, P_1$ are {\it appended to $\Lambda^k$},
and denote this by $P = \Lambda^k*[P_k, P_{k-1}, \dots , P_1]$.
This will be called a [primary] decomposition of $P$.
If $P_k = P_{k-1} = \dots = P_1 = \emptyset$,
then $P$ is a pyramid $\Lambda^k$, and we stop, identifying
$\Lambda^k*[\emptyset, \emptyset, \dots , \emptyset]$ with $\Lambda^k$.
Otherwise we decompose nonempty paths among $P_k, P_{k-1}, \dots , P_1$ in the same way.
Repeating this process recursively, we obtain the complete decomposition of $P$.
Since $P_k, P_{k-1}, \dots , P_1$ are shorter than $P$,
the paths participating in the complete decomposition are
pyramids and empty paths. Thus we call it {\it the complete decomposition of $P$
into pyramids}.

The complete decomposition of a Dyck path can be also represented by a rooted tree:
Given $P = \Lambda^k*[P_k, P_{k-1}, \dots , P_1]$, we represent it
by $\Lambda^k$ as the root with children $P_k, P_{k-1}, \dots , P_1$. Then we
do the same for $P_k, P_{k-1}, \dots, P_1$ and continue recursively, until
all the leaves are pyramids or empty paths. A Dyck paths is easily
restored from its complete decomposition.

\begin{figure}[ht]
$$\resizebox{110mm}{!}{\includegraphics{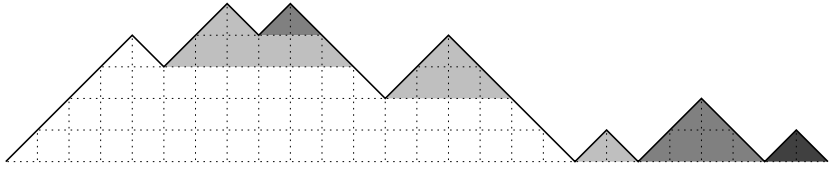}}$$
\caption{Complete decomposition of a Dyck path.}
\label{fig:decomp}
\end{figure}

\begin{figure}[ht]
$$\resizebox{80mm}{!}{\includegraphics{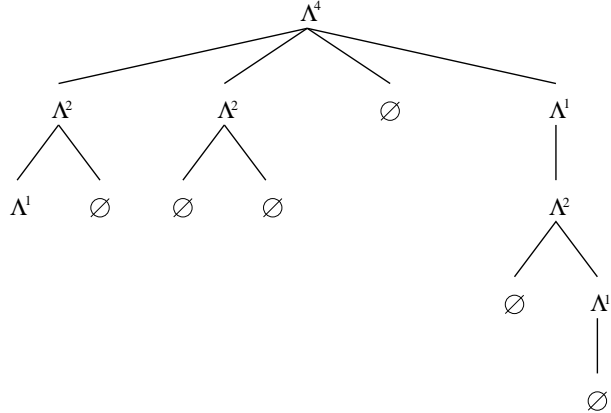}}$$
\caption{Complete decomposition of a Dyck path represented by a rooted tree.}
\label{fig:comptree}
\end{figure}

An example of complete decomposition is
$$\begin{array}{l}
U^4DU^2DUDDDU^2DDDDUDU^2DUD =\\
\qquad=\Lambda^4 * [U^2DUDD, U^2DD, \emptyset, UDU^2DDUD]\\
\qquad=\Lambda^4 * [\Lambda^2 * [UD, \emptyset ], \Lambda^2, \emptyset, \Lambda^1 * [U^2DDUD]] \\
\qquad=\Lambda^4 * [\Lambda^2 * [\Lambda^1, \emptyset ], \Lambda^2, \emptyset, \Lambda^1 * [\Lambda^2*[\emptyset, UD]]]\\
\qquad=\Lambda^4 * [\Lambda^2 * [\Lambda^1, \emptyset ], \Lambda^2,
\emptyset, \Lambda^1 * [\Lambda^2*[\emptyset, \Lambda^1]]],
\end{array}$$
and it is illustrated on Figures~\ref{fig:decomp}
and~\ref{fig:comptree}. On Figure~\ref{fig:decomp} paths appended to
a pyramid are shaded in a more dark colour than the pyramid.

Note that each $U^k$ in a Dyck path results in a $\Lambda^k$ in the complete decomposition.
Therefore the decomposition \eqref{eq:dec} is valid also for coloured Dyck paths:
$\hat{P} = \hat{U}^k D \hat{P}_k D \hat{P}_{k-1} \dots D \hat{P}_1$, and in the complete decomposition
of a coloured Dyck path,
each $U^k$ coloured by $C$ results in $\Lambda^k$ coloured by~$C$.
It is also clear how to restore the coloured Dyck path
from its complete decomposition to coloured pyramids.

We remark that the expression~\eqref{eq:gen2} enumerates thus also
the following structure: rooted trees with $n$ edges, each vertex
$v$ coloured by one of $a_{d(v)}$ colours, where $d(v)$ is the
out-degree of the vertex. For instance, if each vertex $v$ is
coloured by one of $C_{d(v)}$ colours, there are
$\frac{1}{2n+1}\binom{3n}{n}$ such trees.

\subsection {A bijection between non-coloured Dyck paths and NCO-trees}
\label{sec:dycknco}

\begin{figure}[ht]
$$\resizebox{130mm}{!}{\includegraphics{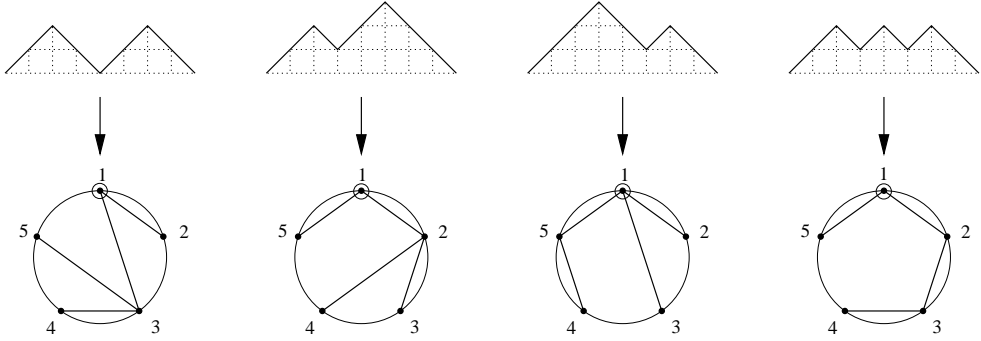}}$$
\caption{The function $\vartheta: \mathcal{D}\rightarrow \mathcal{NCO}$.}
\label{fig:bijdir1}
\end{figure}

We begin with a simple bijection $\vartheta: \mathcal{D} \rightarrow \mathcal{NCO}$.
Given $P \in \mathcal{D}$,
we construct $\vartheta(P)$ according to the following algorithm.
Start with the NC-tree which has one point $1$ and no edges.
Scan $P$ and do the following:
For each $U$, add a new edge beginning in the present point (for a while its end is not determined).
For each $D$, add the next point, move to it and let it be the end of the last incomplete edge.
See Figure~\ref{fig:bijdir1}.

It is easy to see that $\vartheta$ is well defined
(since $P$ is a Dyck path, we never need to complete a non-existing edge,
and all edges are completed by the end),
and that $\vartheta$ is invertible: Given $T \in \mathcal{NCO}$,
scan it from the vertex $1$ clockwise.
Visiting a vertex, first count in-edges incident with it, and then out-edges.
For each in-edge add $U$, for each out-edge add $D$,
and move to the next vertex. It is easy to see that thus obtained Dyck path
$P$ satisfies $\vartheta(P)=T$.
Besides, if $P \in \mathcal{D}(k)$
then $\vartheta(P) \in \mathcal{NCO}(n+1)$.
Thus we have a family of bijections $\vartheta_n: \mathcal{D}(n) \rightarrow \mathcal{NCO}(n+1)$,
for all $n \geq 0$,
which shows in particular that $\vert \mathcal{NCO}(n+1) \vert = C_n$.

\subsection{Definition of $\varphi : \mathcal{D}^{\mathfrak{D}}\rightarrow \mathcal{NC}$}
\label{sec:bij}
First we define $\varphi$ for coloured pyramids.
Consider $\Lambda^k \langle C \rangle$ where $C \in \mathcal{D}$.
We define $\varphi(\Lambda^k \langle C \rangle) = \vartheta(C)$.

\begin{figure}[ht]
$$\resizebox{120mm}{!}{\includegraphics{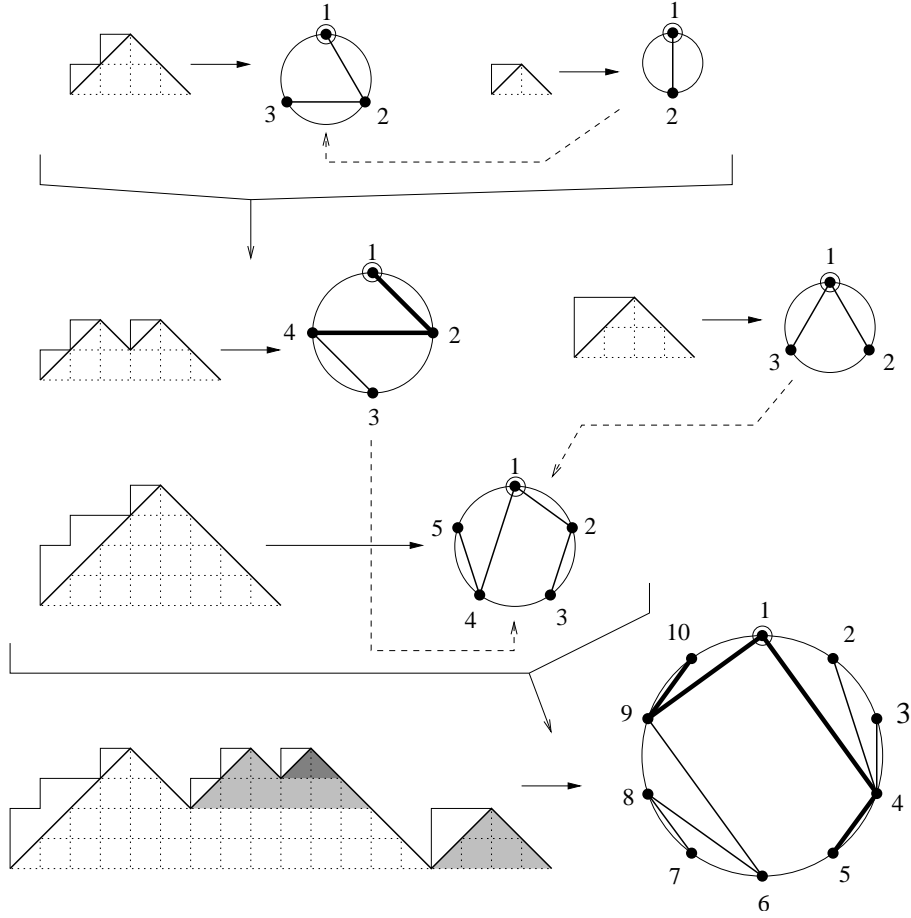}}$$
\caption{The function $\varphi: \mathcal{D}^{\mathfrak{D}}\rightarrow \mathcal{NC}$. The bold edges
are those corresponding to the base pyramid.}
\label{fig:newbij}
\end{figure}

Now we define $\varphi$ for all coloured Dyck paths.
Let $\hat{P}= \hat{\Lambda}^k *
[\hat{P}_{k}, \hat{P}_{k-1}, \dots, \hat{P}_1] \in \mathcal{D}^{\mathfrak{D}}$.
Suppose that we know
$\varphi(\hat{\Lambda}^k)$ and $\varphi(\hat{P}_{i})$ for $i=1, 2, \dots, k$.
For each $i=1, 2, \dots, k$, insert
a copy of $\varphi(\hat{P}_{i})$ into $\varphi(\hat{\Lambda}^k)$
so that
the vertex $1$ of $\varphi(\hat{P}_{i})$
is mapped to the vertex $i+1$ of $\varphi(\hat{\Lambda}^k)$,
and the vertices $2, 3, \dots$ of $\varphi(\hat{P}_{i})$
are mapped clockwise to new vertices between $i$ and $i+1$ in $\varphi(\hat{\Lambda}^k)$
(if $\hat{P}_{i} = \emptyset$ nothing happens).
See Figure~\ref{fig:newbij}.

The function $\varphi$ is invertible.
Given $T\in\mathcal{NC}$,
we want to find $\hat{P} \in \mathcal{D}^{\mathfrak{D}}$
such that $\varphi(\hat{P}) = T$.
Take the subtree of $T$ with root $1$ obtained by recursive adding only out-edges
incident with each reached vertex.
After appropriate relabeling of vertices,
it forms an NCO-tree $V$. It corresponds to the base pyramid $\hat{\Lambda}^k$ of $\hat{P}$,
where $k$ is equal to the number of edges in $V$, and
the colouring is $\vartheta^{-1}(V)$.
For $i = 1, 2, 3, \dots$, the subtree of $T$ attached to the vertex $i$ of $V$
corresponds to $\hat{P}_{i}$ which is determined
recursively. This allows to restore $\hat{P}$.

It is easy to see that if
$\hat{P} \in \mathcal{D}^{\mathfrak{D}}(n)$
then $\varphi(\hat{P}) \in \mathcal{NC}(n+1)$.
Thus we have a family of bijections
$\varphi_n : \mathcal{D}^{\mathfrak{D}}(n)\rightarrow \mathcal{NC}(n+1)$,
for all $n \geq 0$.

This completes the proof of Theorem~\ref{the:dyckdyck}.

\section{Dyck paths coloured by Dyck paths with ascents of bounded length}
\label{sec:mosch}

In this section we restrict the Dyck paths used as colours, considering in this role only Dyck paths
with ascents of bounded length.

Let $\mathcal{M}^m(n)$ be the set of Dyck paths of length $2n$ with
ascents of length $\leq m$. It is known that $\vert \mathcal{M}^m(n)
\vert$ is equal to the $n$-th $m$-generalized Motzkin number. For
$m=2$ we have Motzkin numbers which enumerate Motzkin paths. For $m
\geq n$ we have $\vert \mathcal{M}^{\geq n}(n) \vert= C_n$. In this
sense the sequence of sequences of $m$-generalized Motzkin numbers
``converges'', with $m \rightarrow \infty$, to the sequence of
Catalan numbers.

Among other structures enumerated by $\vert \mathcal{M}^m(n) \vert$ we have
\begin{itemize}
    \item The set of rooted trees on $n+1$ vertices with degree $\leq m$.
    \item The set of all partitions of the vertices of a convex labeled $n$-polygon
    to $(\leq m)$-sets with disjoint convex hulls.
\end{itemize}

Denote $\mathfrak{M}^m= \{\mathcal{M}^m(0), \mathcal{M}^m(1), \mathcal{M}^m(2), \dots \}$.
We consider $\mathcal{D}^{\mathfrak{M}^m}(n)$ for fixed $m$, and prove the following:

\begin{theorem}
\label{the:dyckmo}
There is a bijection between
$\mathcal{D}^{\mathfrak{M}^m}(n)$
and the set of
partitions of the vertices of a labeled convex $(2n)$-polygon
to $(\leq 2m)$-sets of even size with disjoint convex hulls.
The cardinality of both sets is
$\sum_{p=0}^{n/m-1}
\frac{(-1)^p}{n-mp}
\binom{n-mp}{p}
\binom{3n-mp-p}{n-mp-1}$.
\end{theorem}

\subsection{Enumeration}
The generating function $A(x)$ for $\{\vert \mathcal{M}^m(n)\vert \}_{n\geq 0}$
satisfies
\begin{equation}
\label{eq:mo}
A(x) = 1 + xA(x) + x^2A^2(x) + \dots + x^mA^m(x).
\end{equation}

Substituting \eqref{eq:gen2} in \eqref{eq:mo}, we obtain that the generating function
$h_m(x)$ for $\{\vert\mathcal{D}^{\mathfrak{M}^m}(n)\vert\}_{n \geq 0}$ satisfies
\[ h_m(x) = 1 + xh^2_m(x) + x^2h^4_m(x) + \dots + x^mh^{2m}_m(x),\]
which is equivalent to
\[h_m(x) - 1 = xh^3_m(x) - (x h^2_m(x))^{m+1}.\]
Applying the Lagrange inversion formula on
\[h_m(x,a) - 1 = a\left(xh^3_m(x,a) - (x h^2_m(x,a))^{m+1}\right),\]
we obtain
\[ h_m(x,a) - 1 =
\sum_{\ell \geq 1} \frac{a^{\ell}x^{\ell}}{\ell} \sum_{i=0}^{\ell -
1} \sum _{j=0}^{\ell} (-1)^j x^{mj} \binom{3 \ell}{i}
\binom{\ell}{j} \binom{(2m-1)j}{\ell-1-i},\] which implies that the
coefficient of $x^n$ in $h_m(x) = h_m(x,1)$ is
\[[x^n](h_m(x))=\sum_{p=0}^{n/m-1}
\frac{(-1)^p}{n-mp}
\binom{n-mp}{p}
\sum_{i=0}^{n-mp-1}
\binom{3n-3mp}{i}
\binom{2mp-p}{n-mp-1-i}= \]
\begin{equation}
\label{eq:dyckmo}
=
\sum_{p=0}^{n/m-1}
\frac{(-1)^p}{n-mp}
\binom{n-mp}{p}
\binom{3n-mp-p}{n-mp-1}.
\end{equation}

This is the cardinality of $\mathcal{D}^{\mathfrak{M}^m}(n)$.

\subsection{Partitions of convex polygons}
Denote by $\mathcal{E}^m(n)$ the set of
all partitions of the vertices of a convex labeled $(2n)$-polygon
to $(\leq 2m)$-sets of even size with disjoint convex hulls.
We label the vertices of the $(2n)$-polygon by
$[n]\times\{a, b\}=\{a_1, b_1, a_2, b_2, \dots, a_n, b_n\}$
and depict them appearing on a circle clockwise in this order, $a_1$ being the top point.
Denote also $\mathcal{E} = \bigcup_{n\geq 0, m\geq 1}\mathcal{E}^m(n)$.

For all $n\geq 0$, $m\geq 1$ we shall construct a bijection
$\rho_{n, m}: \mathcal{D}^{\mathfrak{M}^m}(n) \rightarrow \mathcal{E}^m(n)$.
It will be presented as a restriction of a bijection
$\rho: \mathcal{D}^{\mathfrak{D}} \rightarrow \mathcal{E}$.

We start with a bijection $\psi$ from $\mathcal{D}$
to $\bar{\mathcal{E}}$, the set of
all partitions of the vertices of a convex labeled polygon
to sets with disjoint convex hulls.
We label the vertices of the polygon by
$1, 2, 3, \dots $ and depict them appearing on a
circle clockwise in this order, $1$ being the top point.

Given $P \in \mathcal{D}$ we construct $\psi(P)$ according to the
following algorithm. Start with the circle without points. Scan $P$
and do the following: For each $U^kD$ add the next point which will
be the first vertex of a $k$-polygon (the other vertices of this
polygon are determined later). For each $D$ not preceded by $U$, add
the next point, move to it and let it be a new vertex of the last
incomplete polygon.

\begin{figure}[ht]
$$\resizebox{130mm}{!}{\includegraphics{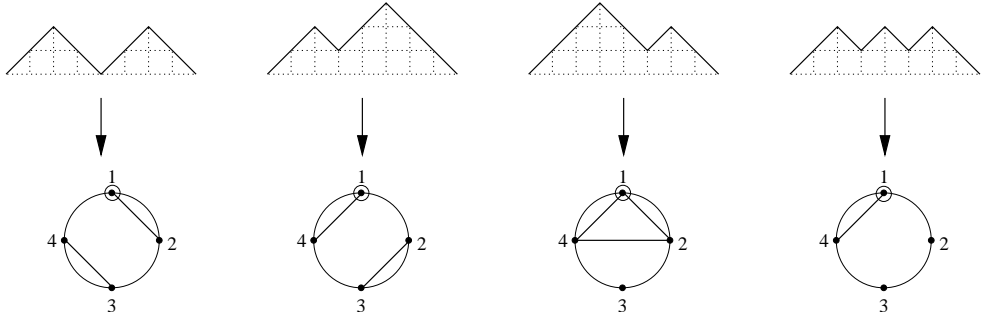}}$$
\caption{The function $\psi: \mathcal{D} \rightarrow \bar{\mathcal{E}}$.}
\label{fig:mo4}
\end{figure}

It is easy to see that $\psi$ is well defined and invertible --
similarly to the function $\vartheta$.
Indeed, comparing Figures~\ref{fig:bijdir1} and~\ref{fig:mo4},
the reader will easily construct a bijection between $\mathcal{NCO}$
and $\bar{\mathcal{E}}$.

Now we define $\rho: \mathcal{D}^{\mathfrak{D}} \rightarrow \mathcal{E}$.
First we define it for coloured pyramids.

Let $M \in \mathcal{D}$.
We define $\rho(\Lambda^k \langle M \rangle)$ to be
a ``duplicated $\psi(M)$'', i.\ e.\
for each polygon with vertices $x_1, x_2, x_3, \dots $ in $\psi(M)$, we have
a polygon with vertices
$a_{x_1}, b_{x_1}, a_{x_2}, b_{x_2}, a_{x_3}, b_{x_3}\dots $ in $\rho(\Lambda^k \langle M \rangle)$.

\begin{figure}[ht]
$$\resizebox{120mm}{!}{\includegraphics{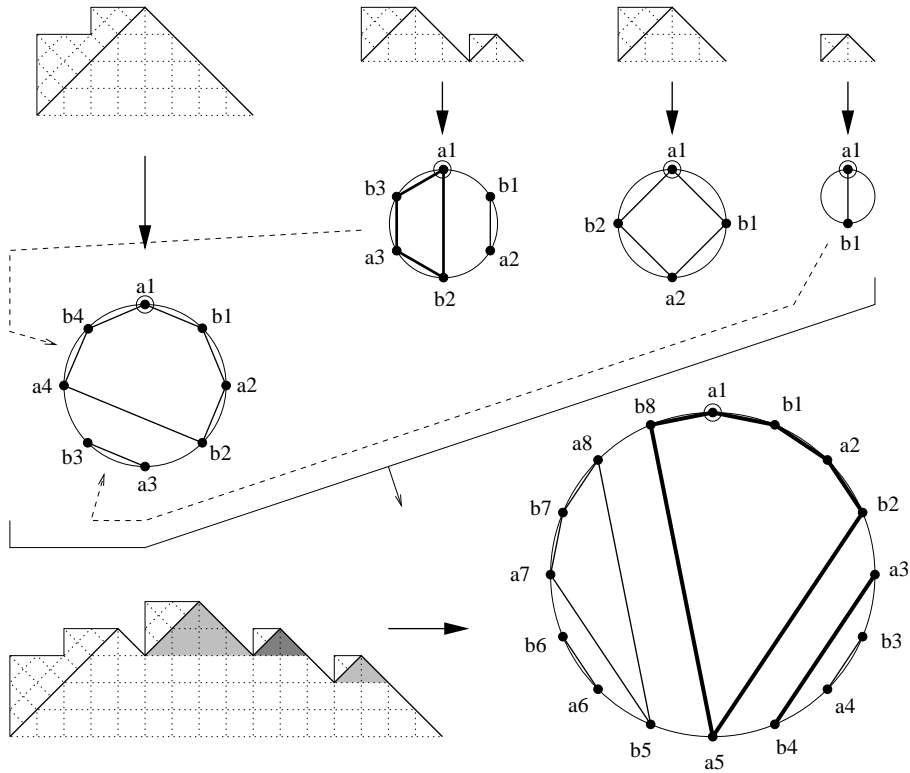}}$$
\caption{The function $\rho: \mathcal{D}^{\mathcal{D}} \rightarrow {\mathcal{E}}$.}
\label{fig:mo3}
\end{figure}

Now let $\hat{P}=\hat{\Lambda}^k * [\hat{P}_k, \hat{P}_{k-1}, \dots, \hat{P}_1] \in \mathcal{D}^{\mathfrak{M}^m}(n)$,
and suppose we know $\rho(\hat{\Lambda}^k)$ and $\rho(\hat{P}_i)$ for $i = 1, 2, \dots, k$.
For each $i = 1, 2, \dots, k$, insert a copy of $\rho(\hat{P}_i)$
into $\rho(\hat{\Lambda}^k)$
so that all the points of $\rho(\hat{P}_i)$ are mapped clockwise to
new points between $a_i$ and $b_i$.
The obtained partition is $\rho(\hat{P})$.
See Figure~\ref{fig:mo3} for an illustration.

The function $\rho$ is invertible.
Given $T\in \mathcal{E}$, we want to find $\hat{P}$
such that $\rho(\hat{P}) = T$.
Consider $T$ as the union of polygons
and choose points of
$T$, beginning from $a_1$ and moving clockwise as follows:
from $a_i$ pass to the vertex connected to it (by an edge of a polygon in the partition),
from $b_i$ pass to $a_{i+1}$.
Denote by $V$
the union of polygons formed by the chosen points after appropriate
relabeling of vertices (these polygons
have bold edges in Figure~\ref{fig:mo3}).
It corresponds to the base pyramid of $\hat{P}$ with colouring
determined by joining points $a_i$ and $b_i$ into one point $i$
and then applying $\psi^{-1}$.
For $i = 1, 2, 3, \dots$, the part of $T$ between the points $a_i$ and $b_i$ of $V$
corresponds to $P_i$ which are determined recursively.
This allows to restore $\hat{P}$.

It is easy to see that if
$\hat{P} \in \mathcal{D}^{\mathfrak{M}^m}(n)$
then $\rho(\hat{P}) \in \mathcal{E}^m(n)$.
In particular, each $k$-ascent in $\hat{P}$
results in a $(2k)$-polygon in $\rho(\hat{P})$.
Thus we have a family of bijections
$\rho_{n, m}: \mathcal{D}^{\mathfrak{M}^m}(n) \rightarrow \mathcal{E}^m(n)$,
for all $n \geq 0, m\geq 1$,
and this completes the proof of Theorem~\ref{the:dyckmo}.

\subsection{Two special cases}
\label{sec:mospec}
We consider two special cases: $m=1$ and $m=n$.

1. Let $m=1$.
Substituting this in~\eqref{eq:dyckmo},
we get
\[
\sum_{p=0}^{n-1}
\frac{(-1)^p}{n-p}
\binom{n-p}{p}
\binom{3n-2p}{n-p-1},\]
which is equal to $C_n$:
since $\vert \mathcal{M}^1(k)\vert = 1$ for each $k$, we have
$\vert \mathcal{D}^{\mathfrak{M}^1}(n) \vert = \vert \mathcal{D}(n) \vert = C_n$.

The partitions of $[2n]$
to even sets with disjoint convex hulls
corresponding to the members of $\mathcal{D}^{\mathfrak{M}^1}(n)$
in the bijection $\rho$ are those in which each set in partition has two members
(all the ways to connect pairs of points of $[2n]$ in convex position by disjoint segments).

2. Let $m=n$.
Substitute this in~\eqref{eq:dyckmo}. The only relevant value of $p$ is $0$,
and we get therefore
\[\frac{1}{n}
\binom{3n}{n-1} =
\frac{1}{2n+1}
\binom{3n}{n},\]
which is expected:
we have $\vert \mathcal{M}^n(n) \vert = C_n$ and thus
$\mathcal{D}^{\mathfrak{M}^n}(n) = \mathcal{D}^{\mathfrak{D}}(n)$.

The corresponding partitions of $[2n]$ are all possible partitions into
even polygons.

\section{Dyck paths coloured by Fibonacci paths}
\label{sec:restr}
In this Section we consider a further restriction of
Dyck paths taken as colours.

Let $\mathcal{F}^{m}(n)$ be the set of Dyck paths of length $2n$
which have the form $\Lambda^{k_1}\Lambda^{k_2}\Lambda^{k_3}\dots$
with $k_i \leq m$ -- a concatenation of pyramids of length no more
than $2m$. Note that $\mathcal{F}^{m}(n) \subset
\mathcal{M}^{m}(n)$. It is known that $\vert \mathcal{F}^{2}(n)
\vert$ is equal to the $(n+1)$-st Fibonacci number; as a
generalization $\vert \mathcal{F}^{m}(n) \vert$ is the $(n+1)$-st
$m$-generalized Fibonacci number (see~\cite[A092921]{S}). The
members of $\mathcal{F}^{m}(n)$ will be therefore called {\it
$m$-generalized Fibonacci paths}. Besides, denote $\mathcal{F}(n) =
\mathcal{F}^{n}(n)$. We have $\vert \mathcal{F}(n) \vert = 2^{n-1}$
for $n>0$, and $\vert \mathcal{F}(0) \vert = 1$.

Denote $\mathfrak{F}^m = \{ \mathcal{F}^m_0, \mathcal{F}^m_1,
\mathcal{F}^m_2, \dots \}$ and $\mathfrak{F} = \{ \mathcal{F}_0,
\mathcal{F}_1, \mathcal{F}_2, \dots \}$. We consider
$\mathcal{D}^{\mathfrak{F}^m}(n)$ for fixed $m$, and
$\mathcal{D}^{\mathfrak{F}}(n)$, and prove the following:

\begin{theorem}
\label{the:dyckfib}
There is a bijection between
$\mathcal{D}^{\mathfrak{F}^m}(n)$and the set of
diagonal dissections of a labeled convex
$(n+2)$-polygon
into $3$-, $4$-, $\dots, (m+2)$-polygons.
The cardinality of both sets is \\
$\sum_{\ell=0}^{n-1}
\frac{1}{\ell+1}
\binom{n+\ell+1}{\ell}
\sum_{i=0}^{\ell+1}
(-1)^i
\binom{n-1-mi}{\ell}
\binom{\ell+1}{i}$.
\end{theorem}

\subsection{Enumeration}
\label{seq:enumfib}
The generating function of the sequence
$\{\vert \mathcal{F}^{m}(n) \vert\}_{n \geq 0}$ is
\[ \sum_{n\geq0}\vert \mathcal{F}(n)^m\vert x^n=\frac{1}{1-x-x^2-\dots-x^m}. \]
Substituting this in \eqref{eq:gen2}, we obtain that the generating
function $g_m(x)$ for the sequence $
\{\vert\mathcal{D}^{\mathfrak{F}^m}(n)\vert \}_{n\geq0}$ satisfies
$$g_m(x)=\frac{1}{1-xg_m(x)-\cdots-(xg_m(x))^m}$$
which is equivalent to
$$g_m(x)-1=x(g_m(x))^2\frac{1-x^m(g_m(x))^m}{1-xg_m(x)}.$$
Applying the Lagrange inversion formula on
$$g_m(x,a)-1=ax(g_m(x,a))^2\frac{1-x^m(g_m(x,a))^m}{1-xg_m(x,a)},$$
we obtain
$$g_m(x,a)-1=
\sum_{\ell\geq0}
\frac{a^{\ell+1}}{\ell+1}
\sum_{j\geq0}
\sum_{i=0}^{\ell+1}
(-1)^ix^{\ell+j+mi+1}
\binom{\ell+j}{j}
\binom{\ell+1}{i}
\binom{2\ell+2+j+mi}{\ell},$$
which implies that the coefficient of $x^n$ in $g_m(x)=g_m(x,1)$ is
\begin{equation}
\label{eq:fib}
\sum_{\ell=0}^{n-1}
\frac{1}{\ell+1}
\binom{n+\ell+1}{\ell}
\sum_{i=0}^{\ell+1}
(-1)^i
\binom{n-1-mi}{\ell}
\binom{\ell+1}{i}.
\end{equation}
This is the cardinality of
$\mathcal{D}^{\mathfrak{F}^m}(n)$.

\subsection{Dissections of a convex polygon}
\label{seq:dissect}
Denote by $\mathcal{R}_m(n)$
the set of all dissections of a labeled convex $n$-polygon
into $i$-polygons with $i=3, 4, \dots, m$ by non-crossing diagonals.
We shall label the vertices of the $n$-polygon by $\alpha, 0, 1, \dots, n$ clockwise,
the top vertex being $\alpha$.
Denote also $\mathcal{R} = \bigcup_{n\geq 0, m\geq 1}\mathcal{R}_{m+2}(n+2)$.

For all $n\geq 0, m\geq 1$
we construct a bijection
$\sigma_{n, m}:\mathcal{D}^{\mathfrak{F}^m}(n)
\rightarrow \mathcal{R}_{m+2}(n+2)$.
It will be presented as a restriction of
a bijection $\sigma:\mathcal{D}^{\mathfrak{F}}
\rightarrow \mathcal{R}$.

\begin{figure}[ht]
$$\resizebox{120mm}{!}{\includegraphics{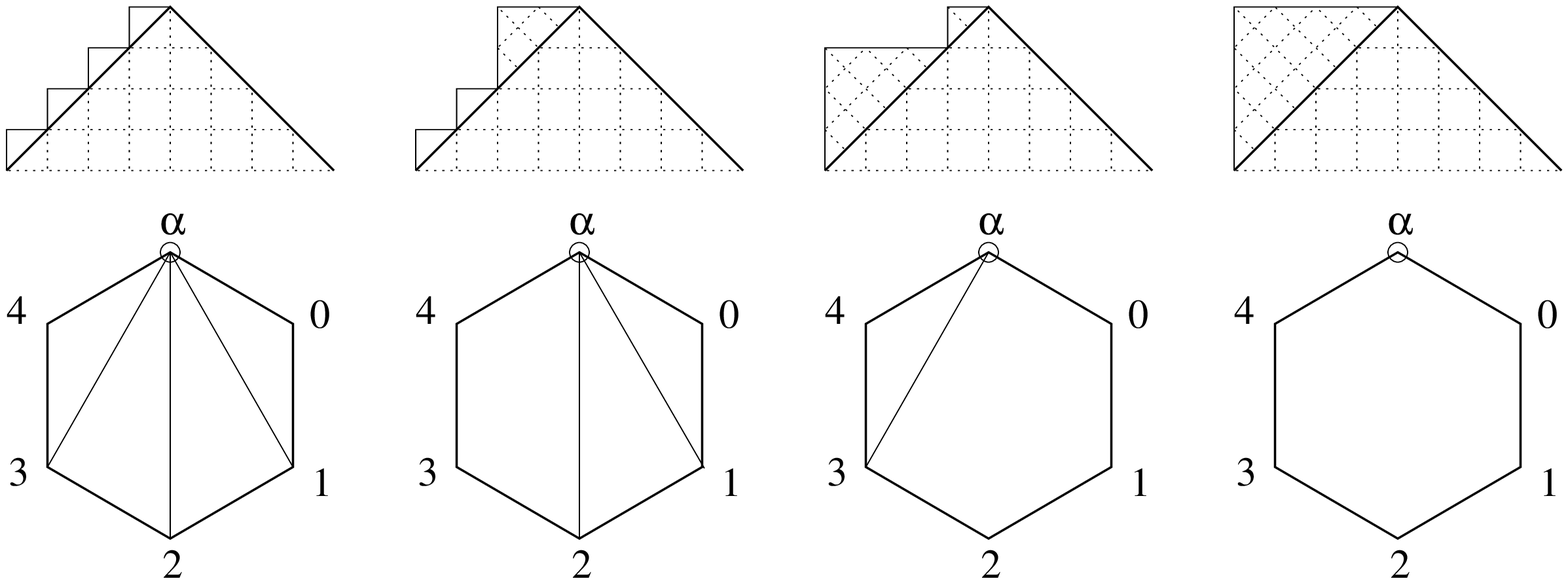}}$$
\caption{Definition of the function $\sigma$ on pyramids.}
\label{fig:pyrdis}
\end{figure}

Consider $\Lambda^k \langle F \rangle$, a pyramid of length $2k$ coloured by
an generalized Fibonacci path $F\in \mathcal{F}(k)$.
We define $\sigma(\Lambda^k \langle F \rangle)$
to be the dissection of the convex polygon with $k+2$ vertices,
taking a diagonal $(\alpha, i)$
if and only if
the path $F$ touches the $x$-axis at point $(i, 0)$
(see Figure~\ref{fig:pyrdis}).

\begin{figure}[ht]
$$\resizebox{100mm}{!}{\includegraphics{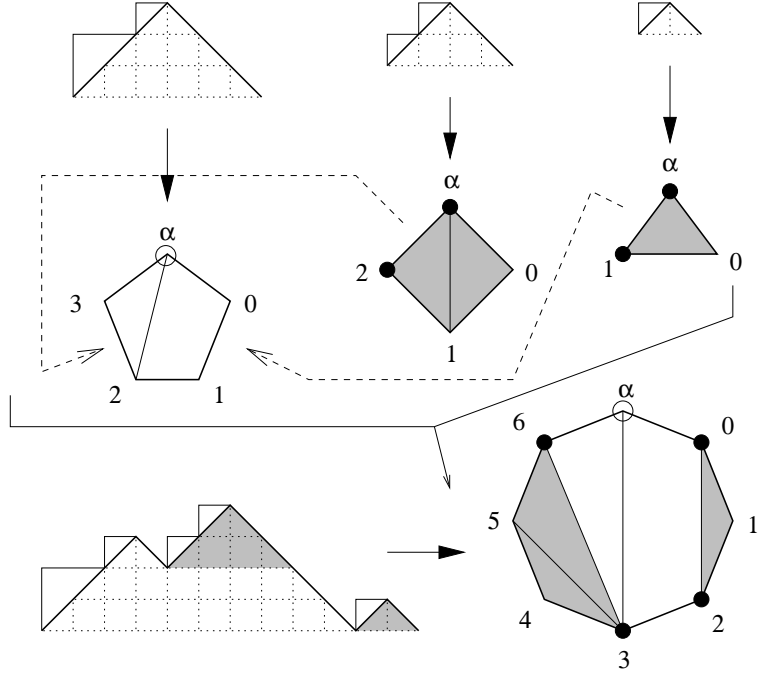}}$$
\caption{The function $\sigma: \mathcal{D}^{\mathfrak{F}}
\rightarrow \mathcal{R}$.}
\label{fig:exdis3}
\end{figure}

Let $\hat{P}=\hat{\Lambda}^k * [\hat{P}_k, \hat{P}_{k-1}, \dots, \hat{P}_1] \in \mathcal{D}^{\mathfrak{F}}$.
Suppose that we know dissections
$\sigma(\hat{\Lambda}^k)$ and $\sigma(\hat{P}_i)$ for $i=1, 2, \dots, k$.
For each $i=1, 2, \dots, k$, insert a copy of $\sigma(\hat{P}_i)$
into $\sigma(\hat{\Lambda}^k)$
so that
the vertex $\alpha$ of $\sigma(\hat{P}_i)$ is mapped to the vertex $i-1$ of $\sigma(\hat{\Lambda}^k)$,
the last vertex of $\sigma(\hat{P}_i)$ is mapped to the vertex $i$ of $\sigma(\hat{\Lambda}^k)$,
and the vertices $1, 2, \dots$ of $\sigma(\hat{P}_i)$ are mapped clockwise to new vertices
between $i-1$ and $i$ of $\sigma(\hat{\Lambda}^k)$.
After relabeling the vertices we
obtain a dissection
$\sigma(\hat{P})$.
See Figure~\ref{fig:exdis3} for an example.

This function $\sigma$ invertible: Let $T\in \mathcal{R}$
and we want to find $\hat{P} \in \mathcal{D}^{\mathfrak{F}}$
such that $\sigma(\hat{P}) = T$.
Let $V$ be the
union of all the polygons in the dissection that have $\alpha$ as a
vertex corresponds, after the appropriate relabeling of its vertices.
$V$ corresponds
to the base pyramid of $\hat{P}$
which is restored immediately.
The part attached to $V$ along the edge $(i-1,
i)$ corresponds to $\hat{P}_i$ which is restored recursively.

Each pyramid of length $2j$ in $F$ (a colouring of $\Lambda^k$)
results in a $(j+2)$-polygon in the dissection of $(k+2)$-polygon.
Therefore if
$\hat{P}\in\mathcal{D}^{\mathfrak{F}^m}(n)$ then
$\sigma(\hat{P})\in \mathcal{R}_{m+2}(n+2)$.
Thus we have a family of bijections
$\sigma_{n, m}:\mathcal{D}^{\mathfrak{F}^m}(n)
\rightarrow \mathcal{R}_{m+2}(n+2)$,
for all $n\geq 0, m\geq 1$,
and this completes the proof of Theorem~\ref{the:dyckfib}.

\subsection{Two special cases} As in Section~\ref{sec:mospec},
we consider two special cases: $m=1$ and $m=n$.
\label{sec:fispec}

1. Let $m=1$.
Substitute this in \eqref{eq:fib}. It can be shown that
$$ \sum_{i=0}^{\ell+1}(-1)^i\binom{n-1-i}{\ell}\binom{\ell+1}{i} = 0$$
for $l\in \{0, 1, \dots n-2\}$. Therefore the whole expression is equal to
$$\frac{1}{n}\binom{2n}{n-1} \sum_{i=0}^{n}(-1)^i\binom{n-1-i}{n-1}\binom{n}{i}=\frac{1}{n}\binom{2n}{n-1}
=C_n,$$
as expected, since
$\mathcal{F}^1(n) = \mathcal{M}^1(n)$ and thus the corresponding
NC-trees are as in Section~\ref{sec:mospec}.
This agrees with a well-known fact that
$\vert \mathcal{R}_3(n+2) \vert$, i.\ e.\
the number of dissections of $n+2$-polygon into triangles, is $C_n$
(see \cite[Page 221, Exercise 19(a)]{S}).

2. Let $m=n$.
Substitute this in \eqref{eq:fib}.
It is clear that
$$(-1)^i
\binom{n-1-ni}{\ell}
\binom{\ell+1}{i}
=0$$
for $i>0$. Therefore the whole expression is equal to
$$\sum_{\ell=0}^{n-1}
\frac{1}{\ell+1}
\binom{n+\ell+1}{\ell}
\binom{n-1}{\ell}
=
\frac{1}{n}
\sum_{\ell=0}^{n-1}
\binom{n+\ell+1}{\ell+1}
\binom{n-1}{\ell}.$$

This expression defines the sequence of ``Little Schr\"{o}der numbers'' \cite[A001003 ]{IS}.
Indeed, it is well known that it enumerates
$\mathcal{R}_n(n+2)$, i.\ e.\
all possible dissections of $n+2$-polygon.
See \cite{callan} for a recent related result.

Let us enumerate $\mathcal{D}^{\mathfrak{F}}(n)$ in another way.
The generating function for $\{\vert \mathcal{F}(n) \vert\}_{n \geq 0}$ is
\[ \sum_{n\geq0}\vert \mathcal{F}(n)\vert x^n=\frac{1-x}{1-2x} = 1 + x + 2x^2 + 4 x^3 + \dots.\]

Substituting this in \eqref{eq:gen2}, we get
$$M = \frac{1-xM}{1-2xM},$$
which is equivalent to
$M-1 = x M (2M - 1)$.
Taking $L=M-1$, we have $L=x(L+1)(2L+1)$,
and by Lagrange's inversion formula,
$$[x^n]L = \frac{1}{n}[L^{n-1}]\left( (L+1)(2L+1) \right)^n=
\frac{1}{n} \sum_{i=0}^{n-1}\binom{n}{i}\binom{n}{i+1}2^{i}.$$
This is the cardinality of $\mathcal{D}^{\mathfrak{F}}(n)$
and thus another expression for Little Schr\"{o}der numbers.

The sequence of Little Schr\"{o}der numbers enumerates ``Little
Schr\"{o}der paths'', see \cite[A001003]{IS}. A {\it Little
Schr\"{o}der path of length $2n$} is a lattice path from $(0,0)$ to
$(2n, 0)$ with moves $U=(1, 1)$, $D=(1, -1)$, $L=(2, 0)$, which does
not pass below the $x$-axis and does not contain an $L$-step on the
$x$-axis.

Denote the set of all Little Schr\"{o}der paths of length $2n$ by $\mathcal{LS}(n)$,
and $\mathcal{LS} = \bigcup_{n \geq 0} \mathcal{LS}(n)$
According to our result, $\vert \mathcal{D}^{\mathfrak{F}}(n) \vert = \vert \mathcal{LS}(n) \vert$.
We construct a simple direct bijection between these sets:

\begin{observation}
\label{the:dycklschr} There is a bijection between
$\mathcal{D}^\mathfrak{F}(n)$ and $\mathcal{LS}(n)$. The cardinality
of both sets is the $n$-th Little Schr\"oder number.
\end{observation}

Let $F \in \mathcal{F}(k)$.
Represent it by a $\{0, 1\}$-sequence $(x_1x_2\dots x_{k-1})$:
consider $F$ as a lattice path from $(0,0)$ to $(2k, 0)$
and let $x_i=1$ if $F$ touches the $x$-axis
at the point $(2i, 0)$, and $x_i=0$ otherwise.

\begin{figure}[ht]
$$\resizebox{130mm}{!}{\includegraphics{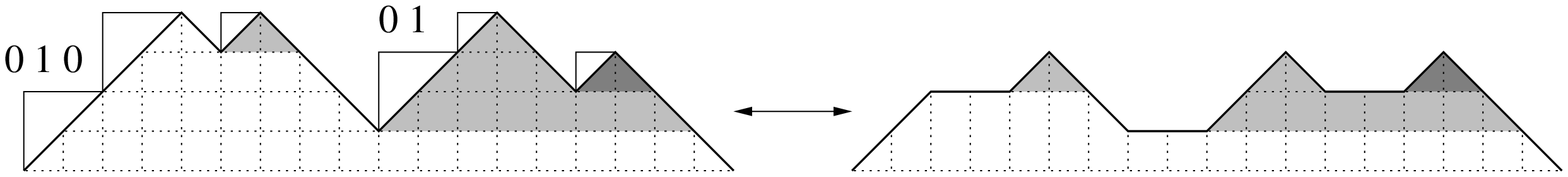}}$$
\caption{The bijection $\mathcal{D}^{\mathfrak{F}}(n) \leftrightarrow \mathcal{LS}(n)$.}
\label{fig:ls2}
\end{figure}

Let $\hat{P} \in \mathcal{D}^{\mathfrak{F}}$.
Consider the complete decomposition of $\hat{P}$.
Each pyramid $\hat{\Lambda}_k$ in this decomposition is
coloured by the members of $\{0, 1\}^{k-1}$.
Replace ${\Lambda}_k \langle (x_1x_2\dots x_{k-1}) \rangle$
with $U^{\beta+1} A_{k-1} A_{k-2} \dots A_2 A_1 D$
where $\beta$ is the number of $1$'s in $(x_1x_2\dots x_{k-1})$,
$A_i=D$ if $x_i=1$, $A_i=L$ if $x_i=0$.
In this way a Little Schr\"oder path is obtained, see Figure~\ref{fig:ls2} for an example.

This function is easily seen to be invertible:
this is based on the fact that any Little Schr\"oder path $P$
may be written in a unique way
as $P= U^\ell X_{k} P_k X_{k-1} P_{k-1} \dots X_2 P_2 D P_1$, where
each $X_i$ is $D$ or $L$, and each $P_i$ is a (possibly empty) Little Schr\"oder path.

Besides, the members of
$\mathcal{D}^{\mathfrak{F}}(n)$ correspond to the members of $\mathcal{LS}(n)$.

\section{Dyck paths coloured by Schr\"oder paths}

Finally, we take $\mathcal{L}_k$ to be the set of all Schr\"oder paths of length $2k$.
A {\it Schr\"oder path of length $2n$} is a
lattice path from $(0,0)$ to $(2n, 0)$
with moves $U=(1, 1)$, $D=(1, -1)$, $L=(2, 0)$,
which does not pass below the $x$-axis.
The set of all Schr\"oder paths of length $2n$ will be denoted by
$\mathcal{S}(n)$.
Schr\"oder sequences are enumerated by \cite[A006318]{IS}.
Let $\mathfrak{S} = \{ \mathcal{S}(0), \mathcal{S}(1), \mathcal{S}(2), \dots \}$.

Denote by $\mathcal{T}(n)$ the set of all lattice paths from $(0,0)$
to $(3n, 0)$ with moves $H=(1, 2)$, $G=(2, 1)$, and $D=(1, -1)$,
that do not pass below the $x$-axis. It is known that
$\mathcal{T}(n)$ is enumerated by \cite[A027307]{IS}. Denote
$t_n=|\mathcal{T}_n|$, for all $n\geq0$.

\begin{observation}
\label{the:dyckschr} There is a bijection between
$\mathcal{D}^\mathfrak{S}(n)$ and $\mathcal{T}(n)$. The cardinality
of both sets is $t_n$.
\end{observation}

We enumerate $\mathcal{D}^{\mathfrak{S}}(n)$ as follows.
The generating function for $\{\vert \mathcal{S}(n)\vert \}_{n \geq 0}$ is
\[ \frac{1-x-\sqrt{1-6x+x^2}}{2x}.\]
Substituting this in~\eqref{eq:gen2} we get
\[M = \frac{1-xM-\sqrt{1-6xM+x^2M^2}}{2xM},\]
or, after simplifications, $M-1 = x(M^3+M^2)$. Denoting $L=M-1$ and
applying Lagrange inversion formula on $L=x(L+1)^2(L+2)$, we finally
get
$$\begin{array}{ll}
[x^n](L)&= \frac{1}{n}[L^{n-1}]((L+1)^2(L+2))^n\\
&=\frac{1}{n}\sum\limits_{i, j \geq 0}\binom{n}{i}\binom{n}{j}
\binom{n}{i+j+1}
2^{i+j+1}=\frac{1}{n}\sum\limits_{k=0}^{n-1}\binom{2n}{k}\binom{n}{k+1}2^{k+1}.
\end{array}$$
This is the cardinality of $\mathcal{D}^\mathfrak{S}(n)$,
and it is known to be an expression for $t_n$, see
\cite[A027307]{IS}.

We show a simple direct bijection between
$\mathcal{D}^{\mathfrak{S}}(n)$ and $\mathcal{T}(n)$.
It can be even said that $\mathcal{T}(n)$ is another representation of
$\mathcal{D}^{\mathfrak{S}}(n)$.

Consider a member of $\mathcal{D}^{\mathfrak{S}}(n)$.
Replace each $k$-ascent by the Schr\"oder path which colours it,
rotated by $45^\circ$ and scaled by $1/\sqrt{2}$
(here we formalize what we already did
in illustrations, see Figure~\ref{fig:dyckdyck}).
We obtain a path with steps
$U=(1,1)$, $N=(0,1)$, $E=(1,0)$, $D=(1, -1)$.

\begin{figure}[ht]
$$\resizebox{130mm}{!}{\includegraphics{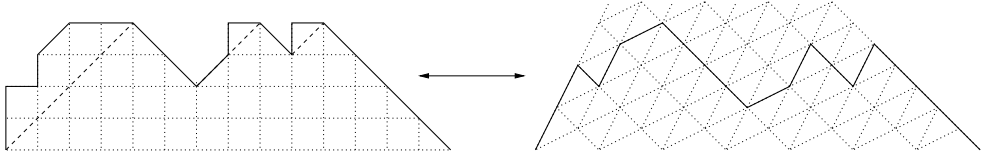}}$$
\caption{The bijection  $\mathcal{D}^{\mathfrak{S}}(n) \leftrightarrow \mathcal{T}(n)$.}
\label{fig:ds2}
\end{figure}

Replacing $N \rightarrow H$,
$U \rightarrow G$,
$E \rightarrow D$,
$D \rightarrow D$, we obtain a member of $\mathcal{T}(n)$.
See Figure~\ref{fig:ds2} for an illustration.

This correspondence is easily seen to be invertible:
in the inverse correspondence, given a member of $\mathcal{T}(n)$, we replace
$H \rightarrow N$,
$G \rightarrow U$,
$D \rightarrow E$ or
$D \rightarrow D$ according to the following rule:
For each $H$, define its {\it match}
to be the closest (from right) $D$
such that the number of $H$'s and of $D$'s between them is equal;
if $D$ is the match of an $H$ then $D \rightarrow E$,
otherwise $D \rightarrow D$.

A restriction of this correspondence is that
between Dyck paths with ascents coloured by Dyck paths
and the members of $\mathcal{T}(n)$
with only moves $H=(1, 2)$ and $D=(1, -1)$.

\end{document}